%% file: multiplicity-paper/multiplicity-paper.tex
\documentclass[12pt]{article}

\usepackage[utf8]{inputenc}
\usepackage[T1]{fontenc}

\usepackage[a4paper, left=2cm, right=2cm, top=1.5cm, bottom=2cm]{geometry}
\usepackage{ragged2e}
\usepackage[pagebackref]{hyperref}
\hypersetup{
colorlinks=true, 
citecolor = blue, 
linkcolor = blue, 
urlcolor = blue, 
anchorcolor = blue
}

\usepackage{amsmath, amsthm, amsfonts}
\usepackage{amssymb, amscd, latexsym}

\usepackage{mathtools} 
\usepackage{caption}

\usepackage{comment, color}
\usepackage{enumitem}

\usepackage{diagbox, array}
\usepackage{tikz, pgf, pgfplots}
\pgfplotsset{compat=newest}

\usepackage{bbm}
\usepackage{lpic}

\usetikzlibrary[patterns]
\usetikzlibrary{arrows, arrows.meta}
\usetikzlibrary{calc}

\definecolor{grey}{rgb}{0.5,0.5,0.5}
\definecolor{forestgreen}{rgb}{0.133,0.545,0.133}
\definecolor{marron}{rgb}{0.6,0.2,0.}
\definecolor{violet}{rgb}{0.5,0.,0.5}
\definecolor{lightpurple}{rgb}{.6,.2,.6}

\theoremstyle{plain}
\newtheorem{theorem}{Theorem}[section]
\newtheorem*{theorem*}{Theorem}

\newtheorem{lemma}[theorem]{Lemma}
\newtheorem{proposition}[theorem]{Proposition}

\theoremstyle{definition}
\newtheorem{definition}[theorem]{Definition}
\newtheorem{example}[theorem]{Example}
\theoremstyle{remark}
\newtheorem{remark}[theorem]{Remark}
\newtheorem{question}[theorem]{Question}
\newtheorem{conjecture}[theorem]{Conjecture}

\newcommand{\bN}{\mathbb{N}}
\newcommand{\bZ}{\mathbb{Z}}
\newcommand{\bQ}{\mathbb{Q}}
\newcommand{\bR}{\mathbb{R}}
\newcommand{\bC}{\mathbb{C}}

\newcommand{\bH}{\mathbb{H}}

\newcommand{\cM}{\mathcal{M}}
\newcommand{\cA}{\mathcal{A}}
\newcommand{\cB}{\mathcal{B}}
\newcommand{\cP}{\mathcal{P}}
\newcommand{\cS}{\mathcal{S}}

\newcommand{\cR}{\mathcal{R}}
\newcommand{\cD}{\mathcal{D}}
\newcommand{\tr}{\mathrm{tr}}
\newcommand{\cW}{\mathcal{W}}

\newcommand{\ch}{\mathfrak{h}}

\newcommand{\mul}{\mathfrak{m}}

\newcommand{\Gap}{\mathrm{Gap}}

\newcommand{\inj}{\mathrm{inj}}

\DeclareMathOperator{\PGL}{PGL}

\DeclareMathOperator{\PSL}{PSL}

\DeclareMathOperator{\vol}{vol}
\DeclareMathOperator{\area}{area}
\DeclareMathOperator{\Isom}{Isom}
\DeclareMathOperator{\axis}{ax}

\title{Large simple length multiplicity on $X(n)$}
\author{Jonah Gaster, D\'idac Mart\'inez-Granado, \\ Christopher-Lloyd Simon}
\date{\today}

\begin{document}

\maketitle

\begin{abstract}
We study multiplicity among lengths of simple closed geodesics on the principal congruence surfaces $X(n)=\Gamma(n)\backslash\mathbb{H}$. We show that for every $\varepsilon>0$ there are positive constants $n_\varepsilon,C_{\varepsilon}$ so that, for $n\ge n_\varepsilon$, there are at least $\exp((\frac{\pi^2}{12}-\varepsilon)n)$ simple closed geodesics on $X(n)$ of the same length $\ell<C_\varepsilon n$. 
We produce the same multiplicity for $\lambda$-length of simple arcs with distinct endpoints 
via the McShane identity, which is immediately relevant as the size of a McShane gap around such an arc is roughly the reciprocal square of its $\lambda$-length.
Our method of proof relies on an independently interesting theorem about the McShane identity for $X(n)$: the total mass of the gaps of pants with two cusps tends to 1 as $n\to \infty$. We present a proof of this latter fact that relies on estimates for the speed of mixing provided by Selberg's 3/16-Theorem (discovered with the help of ChatGPT).

We indicate several avenues for further investigation.

\end{abstract}

\begin{figure}[h]
\centering
\begin{lpic}{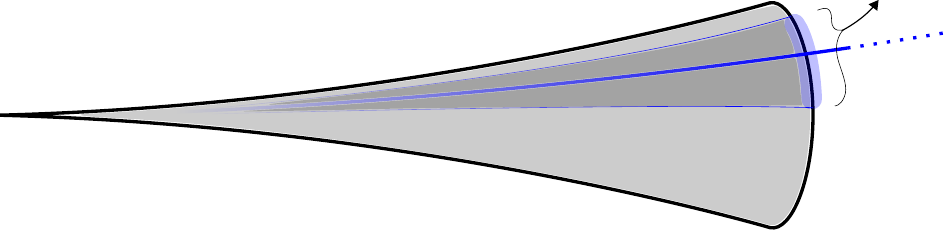(15cm)}
\Large
    \lbl[]{163,42;$\approx1/nq^2$}
\end{lpic}
\caption{A blue simple arc of $\lambda$-length $q$ with distinct endpoints, in a horocycle around the cusp of $X(n)$ of length $n$, contributes a gap of size $\approx1/nq^2$. 
Large multiplicity can be seen as a consequence of the fact that the sum over such arcs of $1/nq^2$ must approximate $n$.}
\label{fig:summary-figures}
\end{figure}

\input{multiplicity-paper/intro}

\input{multiplicity-paper/background}

\input{multiplicity-paper/McShane}

\input{multiplicity-paper/technical-lemma}

\input{multiplicity-paper/proofs}

\newpage
\appendix
\input{multiplicity-paper/Mirzakhani}

\bibliographystyle{alpha} 
\bibliography{biblio.bib}

\end{document}

%% file: multiplicity-paper/intro.tex
\section{Introduction}
In this work, we are concerned with the multiplicity in the simple length spectrum of a hyperbolic surface: how many simple geodesics can have the same length?

Let $Y$ be a connected oriented surface of finite-type with a complete hyperbolic metric. 
A closed geodesic on $Y$ is simple when it is injective (in particular it is primitive).

\begin{definition}
\label{def:simple-length-multiplicity}
Let $\mul(Y,\ell)$ indicate the number of simple closed geodesics on $Y$ of length $\ell$.
The \emph{simple length multiplicity} of $Y$ is given by
\[
\mul(Y)=\sup\{\mul(Y,\ell) \colon \ell 
>0\} \in \bN \cup \{\infty\}~.
\]
\end{definition}

An understanding of $\mul(Y)$ remains almost entirely out of reach.
For instance, the \emph{Markov Unicity Conjecture} is a notoriously difficult open problem which asserts that $\mul(M)=6$ where $M$ is the \emph{modular torus}, 
while the best currently available upper bound is $\mul(M,\ell)=O((\log\ell)^2)$ due to Doan-Li-Nguyen \cite{Doan-Li-Nguyen_simple-length-multiplicity_2026}; note that it is not even known whether $\mul(M)$ is finite!
More generally, one might guess that $\mul(Y)\le 6$ for \emph{every} hyperbolic punctured torus $Y$ \cite{McShane-Parlier_multiplicities_2008}.
Yet more generally, 
the following has been conjectured by Rivin \cite{Schmutz-Schaller_geometry-closed-geodesics_1998}:

\begin{conjecture}[Rivin]
\label{conj:simple-length-mul_bounded-by_chi}
The quantity $\mul(Y)$ is finite and bounded by a function of 
$\chi(Y)$.
\end{conjecture}

We do not touch the thorny question of upper bounds for $\mul(Y)$,
but instead we attempt to construct 
lower bounds that grow quickly with $|\chi|$.
It is trivial to 
construct $Y$ with $\mul(Y)=\Omega(|\chi|)$ using Fenchel-Nielsen coordinates.
The \emph{systole} of $Y$ is the length of
a shortest closed geodesic (necessarily simple), 
and the systole multiplicity is called the \emph{kissing number}. 
One may construct (so called expander) sequences of surfaces $Y$ with large kissing number, demonstrating $\mul(Y)=\Omega(\lvert \chi\rvert^{4/3-\varepsilon})$ \cite{Schmutz-Schaller_geometry-closed-geodesics_1998}.
If one wishes to have dramatically larger constructions, it will be necessary to move beyond systoles: Parlier (and Fanoni-Parlier, in the setting of punctures) has shown that the
kissing number is at most $O(\lvert \chi\rvert^2/\log \lvert \chi \rvert)$ \cite{Parlier_kissing_2013, Fanoni-Parlier_systoles_2015}.

This work concerns the principal congruence surface $X(n)$, namely the quotient of the upper half-plane $\bH$ by the principal congruence subgroup $\Gamma(n)=\ker \: \PSL_2\bZ\to\PSL_2\bZ/n\bZ$.

\begin{theorem}
\label{thm:sup-mul>exp}
For every $\varepsilon>0$ there are $N\in\bN$ and $C_1>1$ so that: for every $n\ge N$ 
there is an $\ell$ with $\frac{1}{C_1}n<\ell<C_1n$ so that
\[
\mul\left( X(n),\ell\right) \ge \exp\left(\left(\tfrac{\pi^2}{12}-\varepsilon\right) n\right)~.
\]
\end{theorem}

It is well-known that $\lvert \chi(X(n))\rvert \asymp n^3$.
Theorem \ref{thm:sup-mul>exp} demonstrates that simple length multiplicity grows faster than any polynomial in $\chi$ along the principal congruence sequence of surfaces, and that any hypothetical bound towards Conjecture~\ref{conj:simple-length-mul_bounded-by_chi} must be 
at least $\exp(\kappa \lvert \chi\rvert^{1/3})$ for some explicit $\kappa>0$.

Our method of proof is nonconstructive: We do not know how to list exponentially many geodesics of the same length on $X(n)$, nor do we know precisely what their length is, although we show that one may find such lengths in $\Theta(n)$. 
Instead we rely on a counting argument involving the McShane identity on $X(n)$ to show existence of many simple geodesic arcs on $X(n)$ with distinct endpoints, and of identical $\lambda$-length. Theorem~\ref{thm:sup-mul>exp} then follows as a result of a well-known relationship between the ($\lambda$-)length of such an arc and the geodesic representative of a small regular neighborhood.

We believe that the same result holds for the expander surfaces coming from congruence subgroups of cocompact arithmetic Fuchsian groups, but our methods should be adapted.

\subsection{Simple arcs and multiplicity in principal congruence surfaces}

Fix $n\in\bN$.
Given $p/q\in\bQ$, we indicate by $\alpha(p/q)$ the projection of the vertical geodesic $\{z\in \bH \colon \Re(z)=p/q\}$ to $X(n)$.
For our purposes, 
the \emph{$\lambda$-length} of $\alpha(p/q)$ is $q$ (geometrically, $\lambda= \exp(\ell/2)$ where $\ell$ is the length of the compact portion of the arc obtained after removing its first and last penetrations in the cusp neighborhoods of area $n$; see \S\ref{sec:background} for more detail).

Observe that the full isometry group of $X(n)$ is $\Isom^\pm(X(n)) \simeq \PGL_2(\bZ)/\Gamma(n)$. 
The next lemma 
follows from the observation that $\Isom^{\pm}(X(n))$ acts transitively on cusps, and the subgroup of $\PGL_2\bZ$ stabilizing $\infty$ is given by $x\mapsto \pm x+k$.

\begin{lemma}
\label{lem:Isom(X(n))-alpha(p/q)}
Each oriented proper geodesic arc on $X(n)$ contains in its $\Isom^\pm(X(n))$-orbit a unique arc of the form $\alpha(p/q)$, oriented from $\infty$ to $p/q$, for a (unique) $p/q\in[0,1/2]$.
\end{lemma}

We are interested in simple arcs and their $\lambda$-length multiplicity, so we propose the following.

\begin{definition}
\label{def:rational spectrum}
The \emph{simple rational spectrum} of $X(n)$ is 
\[
\cS\bQ(n)=\left\{
\frac{p}{q}: \alpha\left(\frac{p}{q}\right) \text{ is simple on } X(n)
\right\}~
\]
and the 
\emph{simple $X(n)$-multiplicity of $q$} is given by
\[
\mul(n,q) = 
\# \left\{
p \in \bN\colon \gcd(p,q)=1,\frac{p}{q}\in\cS\bQ(n) \cap \left[0,\frac{1}{2} \right]
\right\} ~.
\]
\end{definition}

By Lemma~\ref{lem:Isom(X(n))-alpha(p/q)}, $\mul(n,q)$ is equal to the number of $\Isom^\pm(X(n))$-orbits of oriented simple proper geodesic arcs on $X(n)$, of $\lambda$-length $q$, oriented from a fixed cusp.

The simple rational spectrum seems to be largely unexamined. 
We consider 
the present paper as an initial step in this direction.

\begin{example}[$n=2$]
\label{ex:n=2}
The surface $X(2)$ is the unique complete, finite area hyperbolic structure on the three-punctured sphere, and it has only six simple geodesic arcs. More precisely one can check $\cS\bQ(2)\cap [0,1/2] =\{0/1,1/2\}$, so $\mul(2,1)=\mul(2,2)=1$ and $\mul(2,q)=0$ for $q>2$.
\end{example}

For the remainder of the paper, we assume that $n\ge 3$.

\begin{example}[$n=3$]
\label{ex:n=3}
The surface $X(3)$ is conformally isomorphic to $\bC\setminus\{0,1,e^{\pi i/3}\}$, the Riemann sphere punctured at the vertices of the regular tetrahedron.
Since $\Gamma(3)$ is commensurable with the commutator subgroup of $\PSL_2(\bZ)$, the study of $\cS\bQ(3)$ is intimately tied to that of simple proper geodesic arcs on the modular torus, and hence to the Markov numbers \cite{Lehner-Sheingorn_simple-X3-Markov_1984}.

The \emph{Markov numbers} 
\(\cM=\{1,2,5,13,29,\ldots\}~\)
are the integers appearing as solutions to the \emph{Markov cubic} $x^2+y^2+z^2=3xyz$, while their \emph{companions} are defined as $\cM'=\{3m^2:m\in\cM\}$.

The denominators of $\cS\bQ(3)$ consist of $\cM\cup \cM'$, that is $\mul(3,q)>0$ exactly when $q\in\cM\cup \cM'$. 
More precisely, the simple rational spectrum is given by 
\[
\cS\bQ(3) \cap \left[0,\frac 12\right]
=\displaystyle 
\left\{
\frac kn,\,
\frac kn\pm \frac1{3n^2}:\frac kn \text{ is a Markov fraction}
\right\} \cap\left[0,\frac{1}{2}\right]~.
\]
The \emph{Markov Unicity Conjecture} is equivalent to the assertion $\mul(3,q)\le 1$ for all $q\in \cM$.
See \cite{Aigner_Markov-uniqueness-100-years_2013,Gaster_Twist-Numbers-Hyperbolic-Tori_2026,Springborn_worst-approx-ratioals_2024,Veselov_Markov-Fractions-Cohn_2026,Lehner-Sheingorn_simple-X3-Markov_1984,Lehner-Sheingorn_self-intersections-modular-covers_1985,Beardon-Lehner-Sheingorn_Closed-godesics-Markov_1986} for details.
The best upper bound for $\mul(3,q)$ currently available can be deduced from \cite{Doan-Li-Nguyen_simple-length-multiplicity_2026}, which implies that $\mul(3,q)=O((\log\log q)^2)$.
\end{example}

\begin{figure}
\centering
\includegraphics[width=7cm]{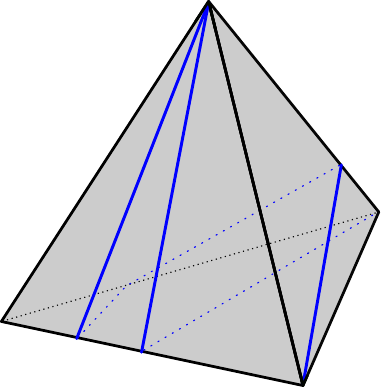}
\hspace{1cm}
\includegraphics[width=5cm]{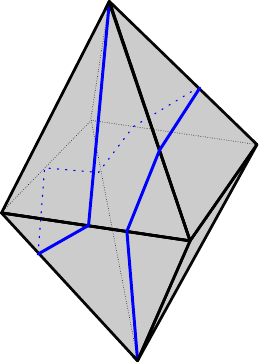}
\caption{On the left: $\frac 12,\frac25\in\cS\bQ(3)$; on the right: $\frac{13}{34}\in\cS\bQ(4)$. 
}
\label{fig:placeholder}
\end{figure}

This seems to be the extent of understanding in the literature for the simple rational spectrum of $X(n)$. 
We plan to initiate a detailed study of $\cS\bQ(n)$ for $n>3$ in a later paper.

\begin{remark}
If the arc $\alpha(p/q)$ is simple on $X(n)$ then it avoids all cusp neighborhoods of area $2$, except the first and last one which it penetrates only at the beginning and end. 
This implies that the continued fraction expansion of $p/q$ has all entries bounded by $n$ 
\cite{Haas_Diophantine-approximation-hyperbolic-surfaces_1986}.
In particular, the study of $\cS\bQ(n)$ refines the question of determining for which integers $q$ is there $p$ such that $p/q$ has continued fraction bounded by $n$, in the spirit of Zaremba's conjecture (see \cite{Zaremba_72}). 
\end{remark}

\subsection{Simple length multiplicity for arcs}

We denote by $\mul_{\ne}(n,q)$ the number of simple arcs on $X(n)$ of the form $\alpha(p/q)$, with $p/q$ reduced and in $[0,1/2]$, and so that the endpoints of $\alpha(p,q)$ are distinct on $X(n)$. 

\begin{remark}
If $n\nmid q$ then the endpoints $\infty$ and $p/q$ are not $\Gamma(n)$-equivalent and so they project to distinct cusps of $X(n)$; in particular $\mul_{\ne}(n,q)=\mul(n,q)$.
\end{remark}

Theorem~\ref{thm:sup-mul>exp} follows immediately from:

\begin{theorem}
\label{thm:m lower bounds}
For every $\varepsilon>0$ there exists $N\in\bN$ and $C_2>0$ so that, for $n\ge N$,
\[
\sup
_{q\le \exp(C_2 n)}
\mul_{\ne}(n,q) \ge \exp\left( \left(\frac{\pi^2}{12}-\varepsilon\right)n\right)~.
\]
\end{theorem}
It can be deduced that the same exact multiplicity statement is true (with a slightly different constant) for multiplicity of simple arcs with both ends in the same cusp: If $\alpha(p/q)$ is simple on $X(n)$ with distinct endpoints, then $\alpha( \tfrac pq \pm \tfrac1{nq^2})$ is simple as well, with same endpoints; the two reduced fractions $\tfrac pq \pm \tfrac1{nq^2}=\tfrac{npq^2\pm1}{nq^2}$ record the two orientations of the same-cusp arc which bounds a once-punctured disk containing $\alpha(p/q)$.

The proof of this theorem relies on the \emph{McShane identity} for a cusped hyperbolic surface.

\begin{theorem*}
[McShane \cite{McShane_identity_1998}]
If $Y$ is a complete finite-area hyperbolic surface with cusp $c$, 
\[
\frac 12 =\sum_{\gamma_1,\gamma_2} \frac1{1+e^{\frac{\ell(\gamma_1)+\ell(\gamma_2)}2}}~,
\]
where $\ell(\gamma_i)$ indicates the hyperbolic length of $\gamma_i$ on $Y$, and the sum is taken over 
unordered pairs $\{\gamma_1,\gamma_2\}$ that cobound an embedded pair of pants with $c$. 
\end{theorem*}

\begin{figure}
    \centering
    \includegraphics[width=5.5cm]{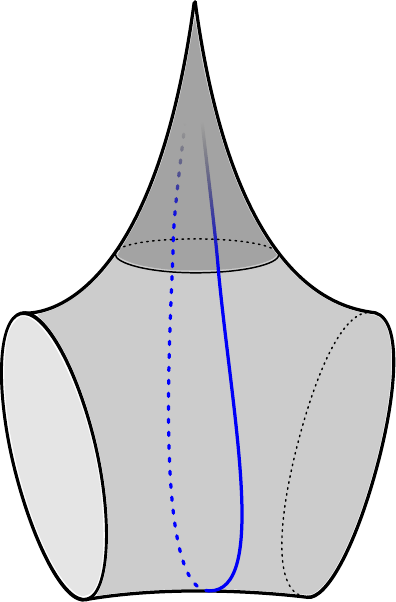}
    \includegraphics[width=8cm]{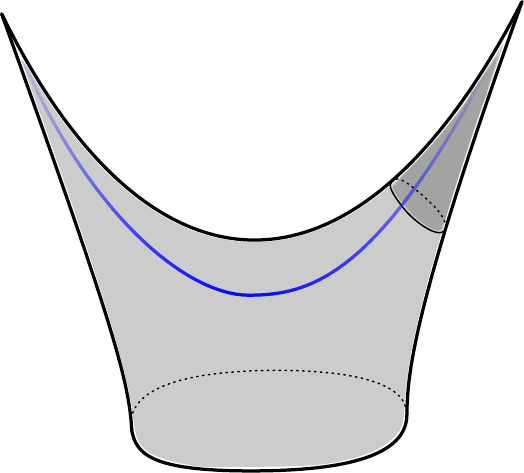}
    \caption{On the left, a pair of pants with one cusp contains a unique simple arc with equal endpoints; the size of the McShane gap at the cusp depends on the $\lambda$-length of the pictured arc, as well as the traces of the two cuff geodesics. On the right, a pair of pants with two cusps contains a unique simple arc with distinct endpoints; the size of the McShane gap at the cusp depends only on the $\lambda$-length of the pictured arc.}
    \label{fig:pants}
\end{figure}

Observe that the McShane identity involves a dichotomy in its summands: Some of the pairs of pants contain one cusp, and others contain two (see Figure~\ref{fig:pants}). 
After multiplying by two,
\[
1=h_=(n)+h_{\ne}(n)~,
\]
where $h_=(n)$ is the portion of the sum coming from pants with one cusp, and $h_{\ne}(n)$ is the portion for two cusps. 
The sums for each of these quantities can be described in precise terms involving the $\lambda$-lengths of certain simple arcs contained in each pair of pants; see Lemma~\ref{lem:McShane Xn}.

One of $h_=(n)$ or $h_{\ne}(n)$ is $\ge 1/2$ infinitely often. 
Let us explain how 
to exploit a lower bound on $h_{\ne}(n)$, and why a lower bound on $h_=(n)$ would not 
be sufficient, deferring more precise statements and justifications to \S\ref{sec:McShane} and \S\ref{sec:proofs}.

On the one hand for $h_{\ne}$, we will find (see \eqref{eq:McShane2}) the estimate
\[
h_{\ne}(n) \approx \sum_{q=1} ^\infty\frac{\mul_{\ne}(n,q)}{n\,q^2}~.
\]

The assumption that $h_{\ne}(n)\ge 1/2$ implies that
\[
\sum_{q=1}^\infty \frac{\mul_{\ne}(n,q)}{q^2} \gtrsim n~.
\]
By definition, for $q\ge 3$ we have $\mul(n,q)\le \phi(q)/2$, where $\phi$ is the Euler totient function. 
By splitting the sum at a large value $Q$, we find
\[
\sum_{q=1}^Q \frac{\phi(q)}{q^2} + \sum_{q>Q} \frac{\mul_{\ne}(n,q)}{q^2} \gtrsim n~.
\]
With the estimate $\sum_1^Q \phi(q)/q^2 = \Theta(\log(Q))$ 
\cite[Ch.~3]{Apostol}
we would nearly arrive at the proof of Theorem~\ref{thm:m lower bounds}: there must exist $C>1$ and $q\gtrsim C^n$ 
so that $\mul_{\ne}(n,q)\gtrsim C^n$. 

On the other hand, 
if one instead 
found that $h_=(n)\ge 1/2$, then 
the expression we obtain for $h_=(n)$ (see \eqref{eq:McShane3}) would yield:
\[
\sum_{\alpha} \frac{\frac1{\tau_1(\alpha)}+\frac1{\tau_2(\alpha)}}{q_\alpha}\gtrsim n~,
\]
where the sum is taken over simple proper geodesic arcs $\alpha$, contained in a pair of pants with \emph{exactly one} (fixed) cusp, and where  $\tau_1(\alpha)$ and $\tau_2(\alpha)$ are the absolute values of the traces of the two geodesic cuffs.
In this case, the numerators depend on more subtle information than merely the $\lambda$-length $q_\alpha$ (namely, as well the traces of the companions of $\alpha$). It is tempting to use the well-known systolic inequality for $X(n)$, which would imply that $1/\tau_i(\alpha)\lesssim 1/n^2$, so that the above estimate becomes
\[
\sum_{q=1}^\infty \frac{\mul(n,q)}{q} \gtrsim n^2~.
\]
One can attempt to estimate as above. However, the partial sums $\sum_q^Q\phi(q)/q=\Theta(Q)$ and $\sum_q^Q\phi(q)/q^2=\Theta(\log Q)$ have dramatically different growth rates. With our methods,  Theorem~\ref{thm:m lower bounds} remains out of reach from the assumption $h_=(n)\ge 1/2 $. 

\begin{question}
Is there a method of reorganizing the sum for $h_=(n)$, so that one finds that the contribution from the terms $\frac1{\tau_1(\alpha)}+\frac1{\tau_2(\alpha)}$ is on average of the order $1/q_\alpha$? (Equivalently, using Lemma~\ref{lem:McShane Xn}, it would suffice to find that $(\tau_1(\alpha)\tau_2(\alpha))^{-1}$ is on the order of $1/nq_\alpha^2$.) 
This could lead to a proof that bypasses the most difficult estimates in the present strategy (in \S\ref{sec:technical}).
\end{question}

\subsection{Controlling single-cusp gaps}

The previous discussion suggests an investigation of the relative sizes of $h_=(n)$ and $h_{\ne}(n)$. 
In order to demonstrate Theorem~\ref{thm:m lower bounds} we use deeper aspects of the geometry of $X(n)$ to show that $ \lim h_=(n)=0$.
Nonetheless, it is reasonable to wonder about the same question, as the hyperbolic structure varies. That is,

\begin{question}
\label{q:relative-sizes}
How do the relative sizes of $h_{=}(Y,c)$ and $h_{\ne}(Y,c)$ vary over $Y$ in the moduli space of complete hyperbolic structures? Same question for $\sum_c h_=(Y,c)$ and $\sum_c h_{\ne}(Y,c)$.
\end{question}
An important observation to make is that nothing uniform in $Y$ can be said: 
If $Y$ has positive genus, one may pinch a convenient pants decomposition, and find a sequence of points in the moduli space so that $\sum_c h_{\ne}(Y,c)$ approaches zero. 

Note that $h_=(Y,c)$ can be analyzed by re-indexing the sum; first, over the set of simple closed geodesics $\cS_Y$, and then over the set of embedded pants $\cP_{=}(c,\gamma)$ with cuff $\gamma\in \cS_Y$ and puncture $c$. 
The set of pairs of pants with cuff $\gamma$ and cusp $c$ is in bijective correspondence with simple orthogeodesics from $\gamma$ to $c$, and the corresponding McShane gap can be expressed as a function of the $\lambda$-length of the corresponding orthogeodesic. 
These considerations are close to Mirzakhani's \emph{bordered} generalization of the McShane identity \cite[Thm.~4.2]{Mirzakhani_Simple-Geodesics-Volumes_2007}. 
For our purposes, the essential consequence of Mirzakhani's version of the McShane identity is:

\begin{proposition}
\label{prop:Basmajian-like sum}
If $Y$ is a connected, oriented, hyperbolic surface of finite-area with $\nu$ cusps $c_1,\dots, c_\nu$, then
\begin{equation*}
    \sum_{i=1}^\nu h_=(Y,c_i) \le  2\sum_{\gamma\in\cS_Y} \frac {\ell(\gamma)}{e^{\ell(\gamma)}-1}~.
\end{equation*}
\end{proposition}

In the case of $X(n)$, whose $\nu(n)$ cusps form a single orbit under $\Isom^{\pm}(X(n))$, Proposition~\ref{prop:Basmajian-like sum} specifies to the claim that $h_=(n)$ is bounded above by $\frac{2}{\nu(n)} \sum_{\gamma\in\cS_n} \frac {\ell(\gamma)}{e^{\ell(\gamma)}-1}$, where $\cS_n$ is the set of simple closed geodesics on $X(n)$.
Since $\nu(n)\asymp n^2$, it is enough to show that 
\[
\sum_{\gamma\in\cS_n} \frac {\ell(\gamma)}{e^{\ell(\gamma)}-1}=o(n^2)
\]
in order to deduce $h_=(n)\to0$.

In fact, toward demonstrating the needed estimate, we prove something slightly stronger.
Recall that on $X(n)$, the set $\cS_n$ of simple geodesics is contained in the set $\cB_n$ of primitive closed geodesics that stay outside horoball neighborhoods of cusps of area $2$.

\begin{proposition}
\label{prop:technical lemma}
Let $B_n$ be the sum given by
\begin{equation}
\label{eq:A_n}    
B_n = \sum_{\gamma\in\cB_n}
\ell(\gamma) e^{-\ell(\gamma)}~.
\end{equation}
We have $B_n=O(n \log n )$.
\end{proposition}

Our proof of Proposition~\ref{prop:technical lemma} relies on a uniform-in-$n$ estimate of exponential mixing for $X(n)$, a consequence of Selberg's celebrated 3/16-Theorem.
We are unsure whether there is a more elementary method of establishing $B_n=o(n^2)$, but would be very interested in obtaining one.
Roughly speaking, the idea of the proof of Proposition~\ref{prop:technical lemma} is to split the sum \eqref{eq:A_n} into three regimes according to whether $\ell(\gamma)$ is on the order of $\log n$ (`short' curves); between $\log n$ and $n \log n$ (`medium' curves); or longer than $n \log n$ (`long' curves). Short curves can be controlled by counting curves on the modular orbifold via class numbers, medium curves are controlled using Selberg's quantified mixing, and the contribution from long curves is shown to be negligible by counting continued fractions with coefficients bounded by $n$.

Finally, we obtain from Propositions~\ref{prop:Basmajian-like sum} and \ref{prop:technical lemma} the desired asymptotic for $h_=(n)$. 

\begin{theorem}
\label{thm:proportion-McShane-gaps}
There is $C_3>0$ so that $h_=(n)\le C_3\log(n)/n$. In particular, $\displaystyle \lim_{n\to\infty} h_=(n)=0$.
\end{theorem}

\begin{remark}
There is a different probabilistic description of the proof of Theorem~\ref{thm:proportion-McShane-gaps} which the reader may find more intuitive, though we've chosen an exposition that passes through the sum $B_n$.
One may interpret the McShane identity as the normalized measure of intervals in a horocyle around $c$ summing to one, where the intervals correspond to outward unit normal vectors whose forward images under the geodesic flow self-intersect \cite{Bridgeman-Tan_identities_2016,Calegari_identities_2010}.

Slightly more precisely, the gaps of pants contributing to either $h_{\ne}(n)$ or $h_=(n)$ measure the subset of points on a horocycle around the cusp where the first self-intersection of an orthogeodesic arc yields either a peripheral or non-peripheral loop, respectively.
The relative sizes of $h_=(n)$ and $h_{\ne}(n)$ are therefore determined by how likely it is that a random orthogeodesic on the horocyclic boundary has first self-intersection that either forms a loop homotopic to a simple closed geodesic, or to a cusp. (There is a small lie here: the event of the outward normal contributing to either $h_=(n)$ or $h_{\ne}(n)$ depends on slightly subtler information, since each pants has two cuffs that must be considered.) 
One can imagine from this description that exponential mixing rates of the geodesic flow are directly relevant, and it turns out that the cusp case dominates as $n\to \infty$.
We do not pursue this direction, but details can be filled in with arguments that are quite close to the control we provide below on $B_n$.
\end{remark}

It would be interesting to have a precise asymptotic for $h_=(n)$.
It is possible to perform a finer dynamical analysis of bounded geodesics of length $\gtrsim n$, using carefully constructed functions on $T^1X(n)$,\footnote{Such functions have been constructed by ChatGPT-Pro.} that in fact demonstrates the stronger bound $h_=(n)=O(1/n)$.
We leave such finer control unaddressed, but consider the following to be a fascinating, subtle question:

\begin{question}
\label{q:growth of h_=}
How does $h_=(n)$ decay to zero? What are good lower bounds for $h_=(n)$?
\end{question}

Returning to Question~\ref{q:relative-sizes}, a potential path towards a more robust version of our main theorem might be to investigate the following:

\begin{question}
Consider an expander family of surfaces $Y(n)$ with given growth assumptions on the genus $g(n)$ and the number of cusps $\nu(n)>0$.
Denote $h_{=}(Y)$ and $h_{\ne}(Y)$ the sum over all cusps $c_i(n)\in Y(n)$ of $h_{=}(Y,c_i)$ and $h_{\ne}(Y, c_i)$: what can we say about $h_{=}(Y(n))/h_{\ne}(Y(n))$? How does it depend on the asymptotic behaviour of $g(n)/\nu(n)$?
\end{question}

Background material and notation are collected in \S\ref{sec:background}, Proposition~\ref{prop:Basmajian-like sum} is proven in \S\ref{sec:McShane}, Proposition~\ref{prop:technical lemma} is proven in \S\ref{sec:technical}, and the passage from Theorem~\ref{thm:proportion-McShane-gaps} to Theorem~\ref{thm:m lower bounds} is in \S\ref{sec:proofs}.

\renewcommand{\contentsname}{Plan of the paper}
\setcounter{tocdepth}{1}

\tableofcontents

\subsection*{Acknowledgements} 

This collaboration was born at the workshop \href{https://www.birs.ca/events/2025/5-day-workshops/25w5411}{Perspectives on Markov Numbers (25w5411)} in Banff International Research Institute, and we thank the institute for the opportunity of organizing and attending this event.

While completing this work,
D\'idac Mart\'inez-Granado was supported by Marie Skłodowska-Curie Action CurrGeo grant (101154865)
and Christopher-Lloyd Simon was supported by the \href{https://perso.univ-rennes1.fr/serge.cantat/ERCGroupsOfAlgebraicTransformations.html}{ERC GOAT 101053021}.

We thank Ben Stucky for participating in the discussions around this project.

\subsection*{Disclosure on the use of AI}

The main ideas (objectives, discoveries and strategies) in this work grew out from the authors in conversations with one another~\cite[Section~4.1(1)]{FGMR2025Markov}, but some of the technical arguments were devised with the interactive help of ChatGPT-Pro. 
In particular, we provided a rough strategy to ChatGPT-Pro who then devised the proof of Theorem~\ref{thm:proportion-McShane-gaps} via the dynamical application of Selberg's theorem, class numbers, and continued fractions.
The authors digested and revised all arguments and take full responsibility for any inaccuracy in this work.

%% file: multiplicity-paper/background.tex
\section{Background on \texorpdfstring{$X(n)$}{X(n)} and \texorpdfstring{$\lambda$-lengths}{lambda-lengths}}
\label{sec:background}

\paragraph{Notation on growth.}

For functions $f,g\colon [0,+\infty) \to [0,+\infty]$, we will write 
\begin{itemize}[noitemsep, align=left]
    \item[$f\lesssim g$] when there is $C>0$ such that $f\le Cg$ (synonymous to Landau's notation $f=O(g)$),
    \item[$f\asymp g$] when $f\lesssim g$ and $g\lesssim f$ (synonymous to Landau's notation $f=\Theta(g)$).
\end{itemize}
Locally defined constants will often be denoted $C$, and will not lose their meaning out of the obvious context. 
The numbered constants $C_i$ will keep their meaning throughout. 

\paragraph{Geometry of $X(n)$.}
We now fix notation and collect well-known information about $X(n)$: 
\begin{enumerate}
    \item The number of cusps of $X(n)$ is $\nu(n)\asymp n^2$.
    \item The volume of $T^1X(n)$ is $V_n\asymp n^3$. 
    \item The length of the systole of $X(n)$ is $s_n=4\log n-O(1)$.
    \item\label{it:horoball packing} The hyperbolic plane $\bH$ admits a $\PSL(2,\bZ)$-invariant packing by horoballs, so $X(n)$ admits a packing by horoball cusp neighborhoods $\{H_j\}$, each of area $n$. 
    \item\label{it:horoball complement} $K_n\subset X(n)$ is the complement of horoball neighborhoods of the cusps, each of area $2$. Every simple closed geodesic of $X(n)$ lies in $K_n$, and for every simple proper geodesic arc $\alpha$, the intersection $\alpha \cap K_n$ is a compact subinterval. 
\end{enumerate}

\paragraph{$\lambda$-lengths.}
We briefly review the notion of $\lambda$-length.

\begin{definition}
\label{def:lambda}
Let $\alpha$ be a proper geodesic arc in $X(n)$ whose ends lie in $H_j$ and $H_k$.
The \emph{$\lambda$-length of $\alpha$, relative to horoball cusp neighborhoods $H_j$ and $H_k$,} is given by 
\[
q_\alpha=\exp\left(\frac{1}{2} \ell(\alpha_0)\right)~,
\]
where $\alpha_0\subset \alpha$ is the segment between the first and last intersections with $\partial H_j$ and $\partial H_k$.
\end{definition}

The following is an elementary calculation, but we have not found this formulation in the literature, so we include a short proof.
An alternate perspective for both identities is that they arise as degenerations of length identities for right-angled hexagons \cite[Chap. 2]{Buser_geometry-spectra-Riemann-surfaces_2010}.
They could also be obtained by derivating Mirzakhani's identity in pairs of pants \cite[Thm.~4.2]{Mirzakhani_Simple-Geodesics-Volumes_2007} (the basic hereditary step in the topological recursion), 
or from identities in \cite{Penner_Decorated-Teichmuller_2012, Fomin-Thurston_cluster-2_2018}.

\begin{lemma}
\label{lem:lambda to length}
Let $\alpha$ be a simple proper geodesic arc on a hyperbolic surface whose ends lie in horoball cusp neighborhoods $H_j$ and $H_k$. 
\begin{enumerate}[label=(\roman*)]
\item If $j\ne k$, then 
\[
\lvert \tr \rvert(\gamma_\alpha) = \area(H_j)\area(H_k) q_\alpha^2-2~,
\]
where $\gamma_\alpha$ is the 
unique simple closed curve that bounds a pair of pants containing $\alpha$.
\item If $j=k$, then
\[
\lvert \tr \rvert(\gamma_{\alpha,1})+\lvert \tr \rvert(\gamma_{\alpha,2})=\area(H_j)q_\alpha~,
\]
where $\gamma_{\alpha,1}$ and $\gamma_{\alpha,2}$ are the two cuffs of the pair of pants containing $\alpha$.
\end{enumerate}
\end{lemma}

\begin{figure}
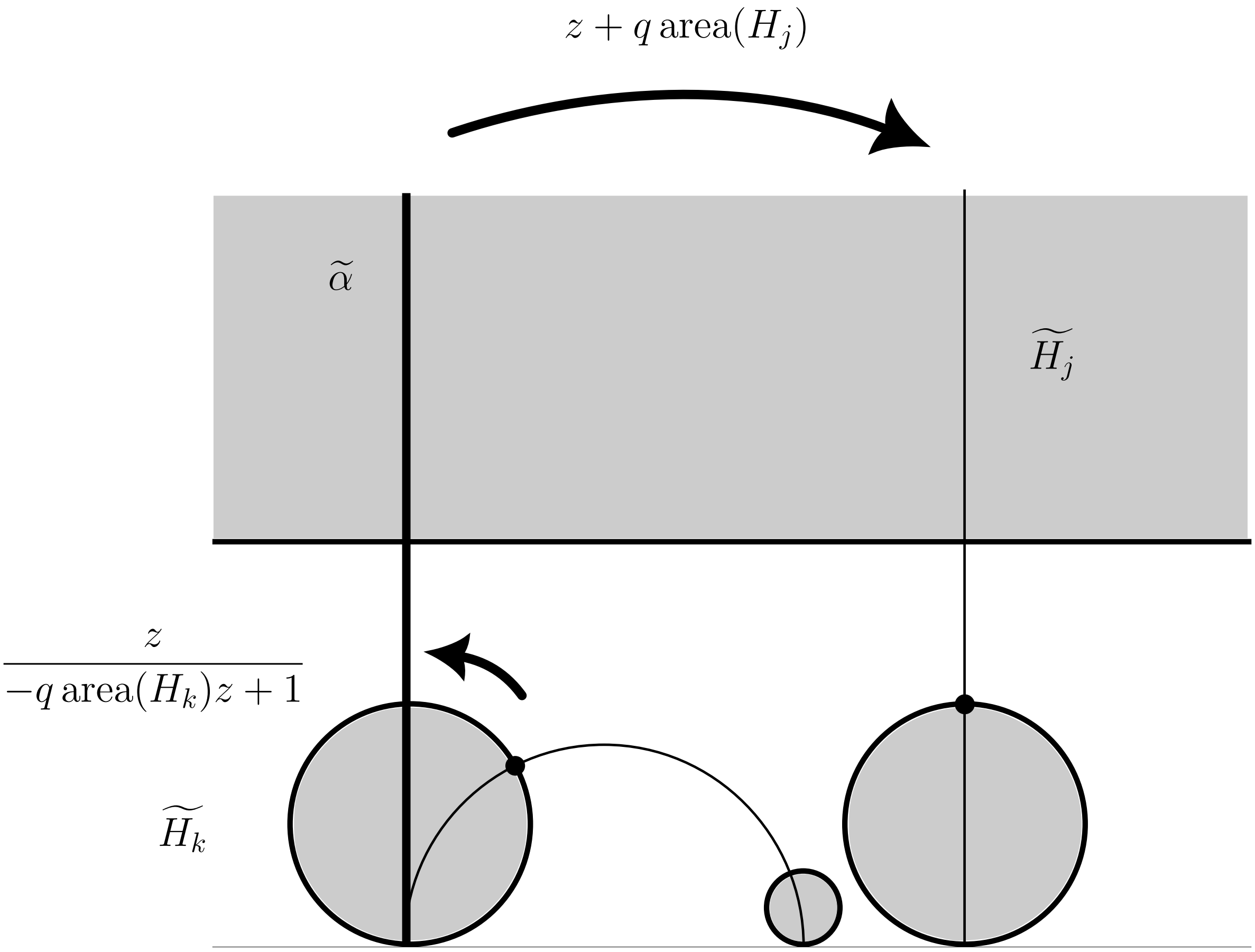

\vspace{1cm}
\centering
\large
\begin{lpic}{images/lambdaTrace1(13cm)}
\lbl[]{48,93;$z+q\area(H_j) $}
\lbl[]{-6,28;$ \displaystyle \frac{z}{-q\area(H_k)z+1}$}
\lbl[]{13,68;$\widetilde{\alpha}$}
\lbl[]{85,60;$\widetilde{H_j}$}
\lbl[]{-3,12;$\widetilde{H_k}$}
\end{lpic}
\caption{If $j\ne k$, the holonomy of $\gamma_\alpha$ is the product of two parabolics.}
\label{fig:lambdaTraces1}
\end{figure}

\begin{figure}
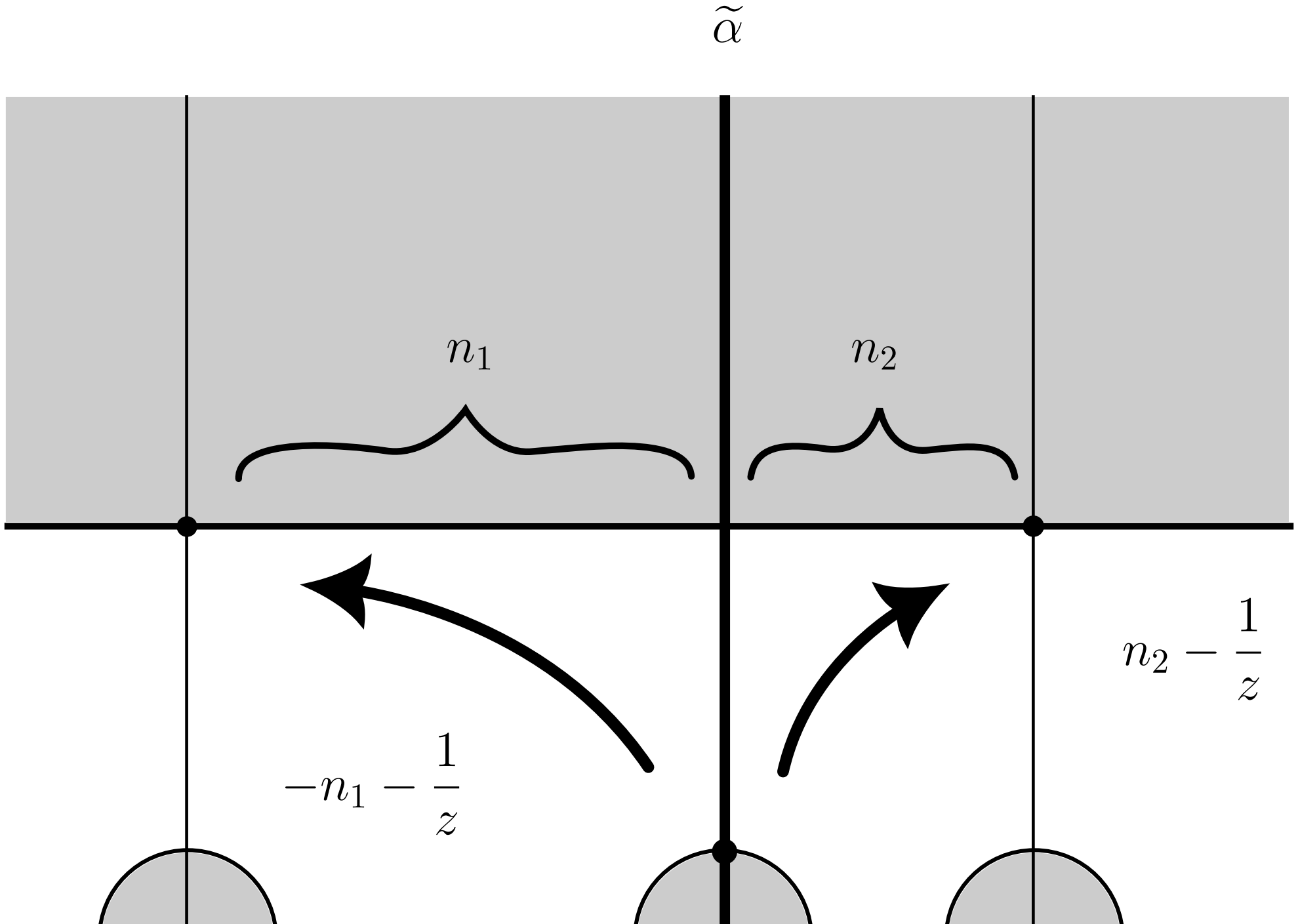

\vspace{1cm}
\centering
\large
\begin{lpic}{images/lambdaTrace2(13cm)}
\lbl[]{30,20;$\displaystyle -n_1 - \frac1z$}
\lbl[]{97,31;$\displaystyle n_2-\frac1z$}
\lbl[]{38,55;$n_1$}
\lbl[]{71,55;$n_2$}
\lbl[]{59,82;$\widetilde{\alpha}$}
\end{lpic}
\caption{If $j=k$, the difference in holonomies of $\gamma_{\alpha,1}$ and $\gamma_{\alpha,2}$ is a parabolic.}
\label{fig:lambdaTraces2}
\end{figure}

\begin{proof}
Let $q=q_\alpha$.
Without loss of generality, we lift $\alpha$ to the imaginary axis in $\bH$, so that $H_j$ and $H_k$ lift to horoballs centered at $\infty$ and $0$, with horocyclic boundary passing through $iq$ and $i/q$, respectively. The holonomies of the oriented simple curves around the cusps at $\infty$ and $0$ are respectively given by (see Figure~\ref{fig:lambdaTraces1})
\[
\begin{pmatrix}1&q\area(H_j)\\0&1\end{pmatrix}
\qquad \text{ and } \qquad
\begin{pmatrix}1&0\\-q\area(H_k)&1\end{pmatrix}~.
\]
Hence, if $j\ne k$ then the holonomy of (an oriented choice for) $\gamma_\alpha$ is conjugate to
\[
\begin{pmatrix}1&q\area(H_j)\\0&1\end{pmatrix}
\begin{pmatrix}1&0\\-q\area(H_k)&1\end{pmatrix}=
\begin{pmatrix}1-q^2\area(H_j)\area(H_k) & q\area(H_j)\\-q\area(H_k)&1\end{pmatrix}~.
\]

If instead $j=k$,
the two (appropriately oriented) cuffs of the pair of pants containing $\alpha$ have holonomies that both take $0$ to $\infty$ and respectively take $i/q$ to $-n_1+qi$ and $n_2+qi$, namely:
\[
\gamma_{\alpha,1} =
\begin{pmatrix}n_1&1\\-1&0\end{pmatrix}
\qquad \text{ and } \qquad
\gamma_{\alpha,2} =
\begin{pmatrix}-n_2&1\\-1&0\end{pmatrix}
\]
and we have $n_1+n_2=q\area(H_j)$, as desired; see Figure~\ref{fig:lambdaTraces2}.
\end{proof}

Throughout this paper, the $\lambda$-length of a proper geodesic arc on $X(n)$ is always taken relative to the maximal cusp neighborhoods $\{H_i\}$ 
fixed above (of area $n$).
For an arc $\alpha$ with distinct endpoints, 
\[
\lvert \tr \rvert(\gamma_\alpha)=n^2 q_\alpha^2-2~.
\]
The passage from Theorem~\ref{thm:m lower bounds} to Theorem~\ref{thm:sup-mul>exp} immediately follows from the one-to-one correspondence $\alpha\leftrightarrow \gamma_\alpha$, between simple arcs $\alpha$ incident to a fixed cusp, with distinct endpoints, and simple closed geodesics $\gamma_\alpha$ bounding a two-cusped pants containing $\alpha$ (see Figure~\ref{fig:pants}): 
Because $\mul_{\ne}(n,q)\le q/2$ for $q\ge 3$,
Theorem~\ref{thm:m lower bounds} yields uniform positive constants $C_4$ and $C_5$ for which there is $q$ between $\exp(C_4n)$ and $\exp(C_5n)$ with  
\[
\mul_{\ne}(n,q) \ge \exp\left( \left( \frac{\pi^2}{12}-\varepsilon\right) n\right)~.
\]
If $\alpha$ is a simple arc with distinct endpoints so that $q_\alpha=q$, then by Lemma~\ref{lem:lambda to length}
\[
\ell(\gamma_\alpha) = 4\log q+4\log n+ O(1)= \Theta(n)~.
\]
Hence
\[
\mul(X(n),\ell) \ge \exp\left( \left( \frac{\pi^2}{12}-\varepsilon\right) n\right)
\]
for some $\ell=\Theta(n)$, as claimed.

%% file: multiplicity-paper/McShane.tex
\section{McShane gaps}
\label{sec:McShane}

\begin{figure}
\centering
\includegraphics[width=6cm]{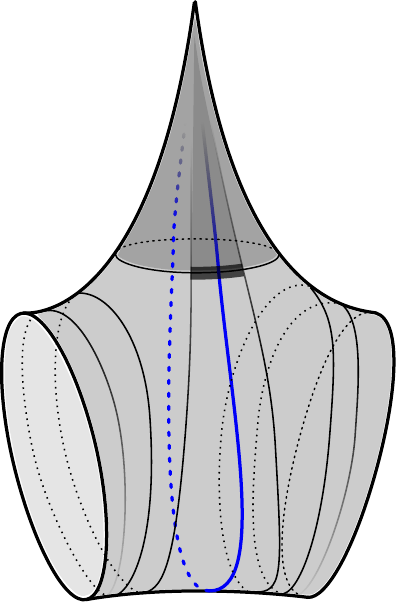}
\includegraphics[width=8cm]{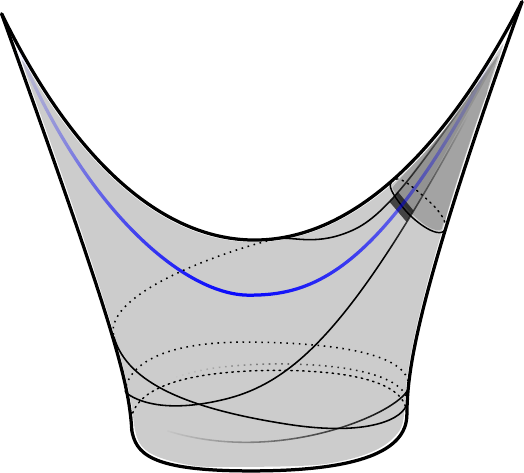}
\caption{The gap $\Gap(P,Y;c)$ for pants $P$ in either $\cP_=(c)$ or $\cP_{\ne}(c)$.}
\label{fig:McShane-gaps}
\end{figure}

Let $c$ be a fixed cusp of $Y$.
Denote by $\cP_=(c)$ the set of embedded pants in $Y$ with a single cusp at $c$ and by $\cP_{\ne}(c)$ the set of embedded pants with two cusps of which one is $c$.
Each pair of pants $P\in\cP_=(c)\sqcup\cP_{\ne}(c)$ has a unique simple geodesic arc $\alpha_P$ with both ends at $c$. We indicate by $\Gap(P,Y;c)$ the size of the corresponding \emph{McShane gap} around $\alpha_P$, given by the normalized length of a certain interval in a horocycle around $c$ in $Y$ demarcating geodesic arcs near $\alpha_P$ that cannot be simple. See Figure~\ref{fig:McShane-gaps} for an illustration, and \cite{Bridgeman-Tan_identities_2016,Tan-Wong-Zhang} for more detail.

Define
\[
h_=(Y,c)= 2\sum_{P\in \cP_=(c)} \Gap(P,Y;c) 
\qquad \text{and} \qquad 
h_{\ne}(Y,c)= 
2\sum_{P\in\cP_{\ne}(c)} \Gap(P,Y;c)
\]
so that we may decompose the McShane identity as
\begin{equation}
\label{eq:McShane0}
1 = h_=(Y,c)+h_{\ne}(Y,c)
\end{equation}

It is possible to give rather concrete expressions for $\Gap(P,Y;c)$ with an important caveat: if $P\in\cP_{\ne}(c)$, $\Gap(P,Y;c)$ depends only on the $\lambda$-length of the corresponding arc, whereas if $P\in \cP_{=}(c)$ then $\Gap(P,Y;c)$ also depends as well on the lengths of the cuffs. This can all be computed quite explicitly on $X(n)$:  

\begin{lemma}
\label{lem:McShane Xn}
Fix a cusp $c$ of $X(n)$.
Let $\cA_{\ne}$ indicate the collection of all simple proper geodesic arcs on $X(n)$ with exactly one endpoint at $c$, and let $\cA_=$ indicate the collection of all simple proper geodesic arcs on $X(n)$ with both endpoints at $c$ that do not bound once-punctured disks on $X(n)$. 
We have
\begin{align}
\label{eq:McShane2}
h_{\ne}(n) &= 2\sum_{\alpha\in\cA_{\ne}} \frac{\omega(n\, q_\alpha)}{n\,q_\alpha} 
\\
\label{eq:McShane3}
h_{=}(n) &=
2\sum_{\alpha\in\cA_=} \frac{\omega(\tau_1(\alpha))+\omega(\tau_2(\alpha))}{n\,q_\alpha}~,
\end{align}
where $\omega(x)=\frac{1}{2}(x-\sqrt{x^2-4})$, $q_\alpha\in\bN$ is the normalized $\lambda$-length of the arc $\alpha$, and, for each $\alpha\in \cA_=$, the integers $\tau_1(\alpha),\tau_2(\alpha)\in\bN_{>2}$ are the two absolute values of traces of the two companion simple closed geodesics of $\alpha\in\cA_=$.
\end{lemma}

\begin{proof}
By Lemma~\ref{lem:lambda to length}, for $\alpha\in\cA_{\ne}$ we have 
\(\lvert \tr\rvert (\gamma_\alpha)=n^2q_\alpha^2-2~,
\) hence the McShane gap is
\[
\Gap(P,Y;c) = 
\frac{1}{1+e^{\frac{\ell(\gamma_1)+\ell(\gamma_2)}2}} = \frac{\omega(nq_\alpha)}{nq_\alpha}~.
\]

By Lemma~\ref{lem:lambda to length}, for $\alpha\in\cA_=$ we have
\(
nq_\alpha = \tau_1(\alpha)+\tau_2(\alpha)\) hence the McShane gap is
\[
\Gap(P,Y;c) = 
\frac{1}{1+e^{\frac{\ell(\gamma_1)+\ell(\gamma_2)}2}} = \frac{\omega(\tau_1(\alpha))+\omega(\tau_2(\alpha))}{nq_\alpha}~. \qedhere
\]
\end{proof}

Now we prove Proposition~\ref{prop:Basmajian-like sum}, recalled below for convenience.

\begin{proposition}
If $Y$ is a connected, oriented, hyperbolic surface of finite-area with $\nu$ cusps $c_1,\dots, c_\nu$, then
\begin{equation*}
    \sum_{i=1}^\nu h_=(Y,c_i) \le  2\sum_{\gamma\in\cS_Y} \frac {\ell(\gamma)}{e^{\ell(\gamma)}-1}~.
\end{equation*}
where $\cS_Y$ is the set of simple closed geodesics on $Y$.
\end{proposition}

\begin{proof}
Denote $g(x,z)=\frac{2}{1+e^{(x+z)/2}}$, and let us unpack the definition of the sum to bound:
\begin{equation*}
    \sum_{i=1}^\nu h_=(Y,c_i) =
    \sum_{i=1}^\nu \sum_{P\in \cP_=(c_i)} 2\Gap(P,Y;c_i)
    = \sum_{i=1}^\nu \sum_{P\in \cP_=(c_i)} 
    g(\ell(\partial^1P), \ell(\partial^2P)).
\end{equation*}
Now re-order the double-sum on the right hand side by summing first over all $\gamma = \partial^1P\in \cS_Y$, then over the cusps $c_i$ and finally over the (nonempty) set of $P\in \cP_{=}(c_i, \gamma)$:
\begin{equation*}
    \sum_{i=1}^\nu h_=(Y,c_i) =\frac12
    \sum_{\gamma\in \cS_Y} h_=(Y,\gamma)
    \quad \mathrm{where} \quad
    h_=(Y,\gamma) =
    \sum_{i=1}^\nu 
    \sum_{P\in \cP_=(c_i, \gamma)} 
    g(\ell(\gamma), \ell(\partial^2P)).
\end{equation*}
(The factor of two arises since, for each $P\in\cP_=(c_i)$ we have $P\in\cP_=(c_i,\partial^1P)\cap \cP_=(c_i,\partial^2P)$.) 
Now fix $\gamma\in\cS_Y$ and let us investigate the corresponding summand $h_=(Y,\gamma)$.

Cut the surface $Y$ along $\gamma$ to obtain a (possibly disconnected) surface $Y(\gamma)$ with $\nu$ cusps and two boundary components of length $\ell(\gamma)$.
We apply \cite[Thm.~4.2]{Mirzakhani_Simple-Geodesics-Volumes_2007} twice, once for each boundary component of $Y(\gamma)$.
Following Mirzakhani's notation, we find that $\ell(\partial Y(\gamma))=2\ell(\gamma)$ can be written as two sums over pants on $Y(\gamma)$, with summands given by Mirzakhani's positive functions $\cR$ and $\cD$. Keeping only the $\cR$ terms for pants incident to a cusp $c_i$ and a component of the boundary $\partial Y(\gamma)$, we find
\begin{equation}
    \sum_{i=1}^\nu 
    \sum_{P\in \cP_=(c_i, \gamma)}
    \cR(\ell(\gamma),0,\ell(\partial^2P))
    \le 2\ell(\gamma)
\end{equation}
where
\[
\cR(x,0,z)=2\log \frac{e^{x/2}+e^{z/2}}{e^{-x/2}+e^{z/2}}~.
\]

Next observe, using $\log u\ge 1-\frac{1}{u}$, that
\begin{equation}
\label{eq:gRBound}
g(x,z)\le \frac{\cR(x,0,z)}{e^x-1}+\frac{\cR(z,0,x)}{e^z-1}~.
\end{equation}
The previous equation indicates that there is a natural splitting of the gap $g(\ell(\gamma),\ell(\eta))$ in two portions, for $\gamma$ and $\eta$ respectively.
From the definition of $h_=(Y,\gamma)$, we find 
\begin{align*}
\frac12\sum_{\gamma\in\cS_Y}h_{=}(Y,\gamma) &= \frac12 \sum_{\gamma\in\cS_Y} \sum_{i=1}^\nu \sum_{P\in\cP_=(c_i,\gamma)} g(\ell(\gamma),\ell(\partial^2P)) \\
&\le \sum_{\gamma\in\cS_Y} \sum_{i=1}^\nu \sum_{P\in\cP_=(c_i,\gamma)}\frac{\cR(\ell(\gamma),0,\ell(\partial^2P))}{e^{\ell(\gamma)}-1} \\
& = \sum_{\gamma\in\cS_Y} \frac1{e^{\ell(\gamma)}-1} \sum_{i=1}^\nu \sum_{P\in\cP_=(c_i,\gamma)} \cR(\ell(\gamma),0,\ell(\partial^2P)) \le \sum_{\gamma\in\cS_Y} \frac{2\ell(\gamma)}{e^{\ell(\gamma)}-1}~. \qedhere
\end{align*}

\end{proof}

%% file: multiplicity-paper/technical-lemma.tex
\section{Controlling the same-cusp gaps}
\label{sec:technical}

In this section we prove Proposition~\ref{prop:technical lemma}, which provides an upper bound for the quantity $B_n$ from \eqref{eq:A_n} which is growing slower than $\nu(n)\asymp n^2$, the number of cusps of $X(n)$.
We decompose $B_n$ into three sums: `short', `medium', and `long' curves. For short curves we use elementary class number bounds and for long curves we use elementary continued fraction bounds. The technical heart of this computation concerns the medium length curves, in which we use uniform exponential-mixing of the geodesic flow on the family of congruence covers $X(n)$.

To begin, 
recall that $\cB_n$ is the set of primitive closed geodesics on $X(n)$ that stay inside $K_n$.
Let $\cB_n(T)$ indicate the elements of $\cB_n$ whose length lies in $[T,T+1)$, and 
\begin{equation*}
    B_n(T)=\sum_{\gamma\in\cB_n(T)} \ell(\gamma) e^{-\ell(\gamma)}~,
    \quad \text{so that} \quad
    B_n=\sum_{T=\lfloor s_n\rfloor}^\infty B_n(T)~,
\end{equation*}
where we recall that the systole of $X(n)$ is $s_n=4\log(n)-O(1)$.

We will bound $B_n(T)$ in three different regimes. 
Let $\kappa >0$ be a constant, whose value will be determined later; see \eqref{eq:K value} for a precise statement. We will say: 
\begin{itemize}
    \item $T$ is \emph{short} provided $T\le 16\log n$; 
    \item $T$ is \emph{medium} provided $ 16\log n<T\le \kappa n\log n$; and 
    \item $T$ is \emph{long} provided $T> \kappa n\log n$.
\end{itemize}

Proposition~\ref{prop:technical lemma} follows immediately from the following bounds:

\begin{lemma}
\label{lem:regimes}
We have: 
\begin{align} 
\sum_{T \text{ is short}} B_n(T) & = O((\log n)^3); \label{eq:B_n1} \\
\sum_{T \text{ is medium}} B_n(T) & = O(n\log n);\label{eq:B_n2} \\
\sum_{T \text{ is long}} B_n(T) & = o(1). \label{eq:B_n3}
\end{align}
\end{lemma}

\subsection{Short curves}
For the short regime, we use elementary class number bounds. Recall the correspondence between conjugacy classes in $\PSL_2(\bZ)$, closed geodesics on the modular surface $X(1)$ and $\PSL_2(\bZ)$-classes of primitive indefinite integral binary quadratic forms \cite{Kontorovich_dynamics-ANT_2016, Sarnak_class-numbers-BQF_1982}.
Denote $\cD$ the set of positive nonsquare discriminants (fundamental or not), and $\ch \colon \cD \to \bN$ the class number.

Conjugacy classes in $\Gamma(n)$ whose translation length is at most $\lceil 16\log n\rceil$ have traces (with sign) confined to $\{2+n^2j\colon j\in\bZ\}$, with $\lvert j\rvert \lesssim n^6$. 
For each element of $\Gamma(n)$ with trace $t$, the discriminant of the corresponding quadratic form is 
$\frac{t^2-4}{u^2}$, where $u$ is the $\gcd$ of the entries of the element, and $n\mid u$. 
Therefore the set of primitive closed geodesics on $X(n)$ of trace $t$ project to at most 
\[
\sum_{n\mid u,\\(t^2-4)/u^2\in\cD} \ch\left( \frac{t^2-4}{u^2}\right)
\]
geodesics on $X(1)$.
The estimate $\ch(D)\lesssim\sqrt D$ implies that the image consists of
$O(|t|\log(2|t|)/n)$ curves on $X(1)$.
Since the cover $X(n)\to X(1)$ has degree $\asymp n^3$, we find at most $O(n^2|t|\log(2|t|))$ geodesics on $X(n)$ of trace $t$. 
(See for instance \cite[Lemma 2.4]{Hashimoto_Selberg-zeta-congruence_2007} for detail.)

As $|t|=2\cosh(\ell/2)$, we have $\ell e^{-\ell}\lesssim \log(2|t|)/t^2$ and, as $T$ is short, $|t|\lesssim n^8$. Therefore 
\[
\sum_{T \text{ is short}} B_n(T) \lesssim n^2 \sum_{t=2+n^2j,\\ 2<|t|\lesssim n^8} \frac{(\log 2|t|)^2}{|t|} \lesssim (\log n)^2 \sum_{1\le j\lesssim n^6} \frac1j~,
\]
and \eqref{eq:B_n1} follows by well-known estimates for partial sums of the harmonic series.

\subsection{Long curves}

Recall that (using the continued fraction expansion encoding of the geodesic flow on the modular surface, see \cite{Series_modular-surface-cont-frac_1985}) a closed geodesic $\gamma$ on $X(1)$ corresponds to a sequence of positive integers $w=(a_1,\ldots, a_{2k})$ of positive even length $2k\in \bN_{>0}$, considered up to even powers of the cyclic shift: we call this a cyclic word and denote it $[w]$.
Under this correspondence, if the closed geodesic $\gamma$ avoids the cusp neighborhood at height $n/2$, then the corresponding cyclic word $[w]$ has all its entries $a_i\le n$.

A closed geodesic $\gamma\in\cB_n$ on $X(n)$ projects to such a closed geodesic in $X(1)$, and hence determines such a cyclic word $[w]$. (Note that the primitive curve $\gamma$ may result in a non-primitive word $w$; in $\gamma\mapsto[w]$ we retain the full period for the word $w$ to ensure the map from curves to words preserves length.)
Since the cover $X(n)\to X(1)$ has degree $\asymp n^3$, there are $O(n^3)$ geodesics in $\cB_n$ yielding the same cyclic words.
The strategy below to obtain \eqref{eq:B_n3} is to bound the contributions in $B_n$ that arise from $[w]$.

Note that any word $w=(a_1,\dots,a_k)$ with entries $a_i\in \bN_{>0}$ yields a matrix $M_w \in\PGL_2(\bZ)$: 
\[
M_a=\begin{pmatrix}a&1\\1&0\end{pmatrix}, \ \ \ \text{ and }
\qquad
M_w=M_{a_1}\cdots M_{a_k}~.
\]
We declare the \emph{length} of the word to be given by $\ell(w)=2\log \rho(M_w)$, where $\rho$ is the spectral radius; when $w$ is the even word associated to the curve $\gamma$, we have $\ell(w)=\ell(\gamma)$. 

Let $I_w\subset(0,1)$ indicate the subset whose continued fraction expansions begin with $w$. 
Each digit $a_i$ in $w$ narrows the size $\lvert I_w\rvert$; the larger the digit, the narrower the interval.
The following are elementary; we provide a proof sketch.

\begin{lemma}
\label{lem:continued fraction computation}
Let $w=(a_1,\ldots,a_{k})$.
We have
\begin{enumerate}[label=(\roman*)]
    \item $e^{-\ell(w)}\le 2\lvert I_w \rvert$ and
    \item $\displaystyle \ell(w)\le 2\sum_i \log(1+a_i)~.$
\end{enumerate}
\end{lemma}

\begin{proof}
Given word $w=(a_1,\ldots,a_k)$, let $q_i$ indicate the $i$th denominator of the convergents. 
Observe that 
\[
|I_w|=\frac1{q_k(q_k+q_{k-1})}~.
\]
Moreover,
\[
M_w=\begin{pmatrix}q_k&q_{k-1}\\p_k&p_{k-1}\end{pmatrix}, \ \ \ \text{ and }
\qquad
\rho(M_w)\ge q_k~.
\]
Therefore $e^{-\ell(w)}\le q_k^{-2}\le 2|I_w|$.
Moreover, 
\[
\rho(M_w)\le \|M_w\|_\infty \le \prod_i \|M_{a_i}\|_\infty = \prod_i(1+a_i)~.
\]
\end{proof}

One may show better: the following implies that the set of numbers whose continued fraction is bounded by $n$ has Hausdorff-dimension strictly less than $1$ (as is well known \cite{Shallit_Reals-bounded-partial-quotients_1992, Hensley_Continued-fraction-Cantor-sets_1992}).

\begin{lemma}
\label{lem:continued fraction harder computation}
There exists an absolute constant $\beta>0$ so that for all $n$, $\displaystyle \eta_n=\frac \beta{n+1}$ satisfies 
\[
\sum_{w\in\{1,\ldots,n\}^k} |I_w|\prod_{i=1}^k (1+a_i)^{2\eta_n}\le \left( \frac {2n+1}{2n+2}\right)^k~.
\]
\end{lemma}

\begin{proof}
Let us indicate the word obtained by appending digit $a$ to $w$ by $wa=(a_1,\ldots,a_k,a)$.
The relative weight of appending $a$ is given by
\begin{equation}
\label{eq:intermediate_lem33}
\frac{|I_{wa}|}{|I_w|} = \frac{1+r}{(a+r)(1+a+r)}~\leq \frac{2}{a^2},
\end{equation}
where $r=q_{k-1}/q_k\in[0,1]$.

We may thus find the relative weight of the regions excluded from $|I_w|$, so that 
\[
\sum_{a=1}^n \frac{|I_{wa}|}{|I_w|}\le 1-\frac1{n+1}~.
\]

We may also estimate a relative mass obtained by appending $1\le a\le n$, weighted by a power of $(1+a)$, with $\eta$ to be chosen momentarily:
\begin{align*}
\sum_{a=1}^n \frac{|I_{wa}|}{|I_w|}(1+a)^{2\eta} & = \sum_{a=1}^n \frac{|I_{wa}|}{|I_w|} \ + \ \sum_{a=1}^n \frac{|I_{wa}|}{|I_w|}\left((1+a)^{2\eta}-1\right) \\
& \le 1-\frac1{n+1} + \sum_{a=1}^\infty \frac{|I_{wa}|}{|I_w|}\left((1+a)^{2\eta}-1\right)~.
\end{align*}
Observe that, provided $0\le \eta\le \frac14$ we have:
\[
(1+a)^{2\eta}-1 = \int_0^\eta 2\log(1+a)(1+a)^{2s}ds \leq 2\eta\log(1+a)\sqrt{1+a}~.
\]
Therefore~\eqref{eq:intermediate_lem33} implies that there is a constant $C$ such that
\[
\sum_{a=1}^\infty \frac{|I_{wa}|}{|I_w|}\left((1+a)^{2\eta}-1\right) \le \sum_{a=1}^\infty \frac2{a^2}\left((1+a)^{2\eta}-1\right) \le C\eta~.
\]
Therefore
\[
\sum_{a=1}^n \frac{|I_{wa}|}{|I_w|}(1+a)^{2\eta} \le 1-\frac{1}{n+1}+C\eta~.
\]
We choose a sufficiently small absolute constant $\beta>0$ so that after setting $\eta=\beta/(n+1)$ we have $C\eta\le \frac1{2(n+1)}$. Hence
\[
\sum_{a=1}^n \frac{|I_{wa}|}{|I_w|}(1+a)^{2\eta} \le 1-\frac1{2(n+1)}~.
\]
Iterating this estimate $k$ times, starting with the full interval $(0,1)$, proves the claim.
\end{proof}

Now Lemma~\ref{lem:continued fraction computation} implies that
\[
e^{-(1-\eta_n)\ell(w)}=e^{-\ell(w)}e^{\eta_n\ell(w)}\le 2|I_w|\prod_i (1+a_i)^{2\eta_n}~,
\]
so that Lemma~\ref{lem:continued fraction harder computation} then implies that
\[
\sum_w e^{-(1-\eta_n)\ell(w)} \le 2\sum_{k\ge 1} \left(\frac{2n+1}{2n+2}\right)^k\lesssim n~.
\]

We connect this to the long regime of $B_n$. For each $\ell>0$,
\[
\ell e^{-\ell}=\left(\ell e^{-\eta_n\ell/2}\right) e^{-(1-\eta_n)\ell} e^{-\eta_n\ell/2} \le \frac2{e\eta_n} e^{-(1-\eta_n)\ell} e^{-\eta_n\ell/2}~,
\]
so that
\[
\sum_{\ell(w)\ge L} \ell(w)e^{-\ell(w)} \lesssim \frac1{\eta_n}e^{-\eta_n L/2} \sum_w e^{-(1-\eta_n)\ell(\gamma)}\lesssim n^2 e^{-\beta'L/n}~,
\]
for some absolute constant $\beta'>0$.

Because there are $\lesssim n^3$ geodesics in $\cB_n$ for a given cyclic word $w$, we find that
\[
\sum_{T=L}^\infty B_n(T) \lesssim n^3 \sum_{\ell(w)\ge L} \ell(w)e^{-\ell(w)} \lesssim n^5 e^{-\beta' L/n}~.
\]

The sum for long curves is obtained by setting $L=\kappa n \log n$. Observe that 
\[
5\log n - \frac{\beta' L}{n} = 5\log n - \frac{\beta'\kappa n \log n}{n}=(5-\beta'\kappa)\log n~.
\]
Hence, provided we choose $\kappa$ so that 
\begin{equation}
\label{eq:K value}
\beta'\kappa >5~,
\end{equation}
the result \eqref{eq:B_n3} follows.

\subsection{Medium curves}
We turn to medium length curves. Throughout this section, let $X=X(n)$ and projection $\pi:T^1X\to X$ the natural projection. We indicate the normalized Liouville measure of $T^1X$ by $\mu$, and the unnormalized area measure on $X$ by $A$.

The quantity $B_n(T)$ admits a pleasant geometric formulation: $\ell(\gamma)e^{-\ell(\gamma)}$ is a rough estimate for the volume of a certain tubular neighborhood in $T^1X$ of a simple closed geodesic, whose projection to $X$ forms an embedded collar neighborhood (cf.~the Collar Lemma~\cite{Keen_Collars-Riemann-surfaces_1974}).
In fact, the regions we will need are not precisely open neighborhoods projecting to the annuli from the Collar Lemma (indeed, curves in $\cB_n$ may even be not simple, in which case the relevant reference is \cite{Basmajian}), but it is nonetheless helpful to think of them as such.
We isolate the precise statement that is needed in our control of $B_n(T)$ below.

\begin{lemma}
\label{lem:tubular nbhd}
There are constants $C, r, \delta>0$ such that: for every $T$ there are tubular neighborhoods $R_\gamma\subset T^1X$ indexed by $\gamma\in\cB_n$, so that the following hold:
\begin{enumerate}[label=\roman*)]
\item\label{it:inj}
For all $v\in R_\gamma$ we have $\inj_{X}(\pi v)\ge r$.
\item\label{it:thick} 
For every $v\in R_\gamma$ and $t>0$ such that $|t-\ell(\gamma)|\le \delta$ we have \(d_X(\pi v,\pi g^tv)<r/2\).
\item\label{it:vol} 
For every $\gamma\in\cB_n$ we have
$\vol R_\gamma \ge C \ell(\gamma) e^{-\ell(\gamma)}$.
\item\label{it:disjoint} 
For distinct $\gamma,\eta\in\cB_n(T)$ we have
$R_\gamma\cap R_\eta=\emptyset$.
\end{enumerate}
\end{lemma}
\begin{proof}
Recall that $K_n \subset X$ is the complement of the open cusp neighborhoods of area 2. 
For all $n$, the injectivity radius at every point in $K_n$ is $\geq \sinh^{-1}(1)$.

Fix $\rho,\theta>0$.
For each $\gamma\in \cB_n$ of length $\ell$, its holonomy $M_\gamma\in\PSL_2(\bR)$ has geodesic axis $\axis(M_\gamma)$, and we form $\widetilde{R_\gamma}$ by taking the interval of directions at points $z$ in the $\rho$-neighborhood of $\axis(M_\gamma)$, centered on the direction from $z$ to $M_\gamma(z)$, of angular width $2\theta e^{-\ell}$. Let $R_\gamma\subset T^1X$ be the projection of $\widetilde{R_\gamma}$ under the covering map.
Note that it is immediate that, for every $v\in R_\gamma$, the projection $\pi v$ lies in the $(\sinh^{-1}1-2\rho)$-thick part of the surface; because $\sinh^{-1}1>1/2$ we may choose $\rho$ sufficiently small so that the first item is satisfied for  $r\in (1/2, \sinh^{-1}(1))$.

Changes in the initial direction of size at most $\theta e^{-\ell}$ changes the endpoint, after flowing for time $\ell$, by $O(\theta)$, uniformly in $\ell$. 
Thus every vector in $R_\gamma$ returns within $O(\rho+\theta+\delta)$ of its starting point for $|t-\ell|\le \delta$, and the second item is fulfilled for $\theta$, $\rho$, $\delta$ sufficiently small. 

Suppose $v\in R_\gamma\cap R_\eta$, choose a lift $\widetilde{v}\in \widetilde{R_\gamma}\cap \widetilde{R_\eta}$, and let $z=\pi \widetilde{v}$. 
Let $t_\gamma=d_\bH(z,M_\gamma(z))$ and $t_\eta=d_\bH(z,M_\eta(z))$, each within $2\rho$ of $\ell(\gamma)$ and $\ell(\eta)$, respectively, by the triangle inequality. We find that
\[
d_\bH(M_\gamma(z),\pi g^{t_\gamma}v) = O(\rho+\theta+\delta) \quad \text{ and } \quad d_\bH(M_\eta(z),\pi g^{t_\eta}v) = O(\rho+\theta+\delta)~.
\]
Because
$|t_\gamma-t_\eta|<1+2\rho$, $M_\gamma(z)$ and $M_\eta(z)$ are within
$1+O(\rho+\theta+\delta)$
of each other, so
\[
d_\bH(z,M_\gamma^{-1}M_\eta z) \le 1 + O(\rho+\theta+\delta)~.
\] 
If $M_\gamma\ne M_\eta$, then 
\[
d_\bH(z,M_\gamma^{-1}M_\eta z)\ge 2\inj_X(\pi z)\ge 2r >1~.
\]
Hence for $\theta$, $\rho$, and $\delta$ sufficiently small we find that $M_\gamma=M_\eta$, and item four is satisfied. 

Note that the same argument demonstrates that $\widetilde{R_\gamma}/\langle M_\gamma\rangle \to R_\gamma $ is injective provided $\gamma$ is primitive, since any covering transformation $g$ with $g\cdot \widetilde{R_\gamma}\cap \widetilde{R_\gamma}\ne\emptyset$ must centralize $M_\gamma$, and hence lie in $\langle M_\gamma\rangle$.
Therefore it is straightforward to compute $\vol(R_\gamma)$:
\[
\vol(R_\gamma)=2\theta e^{-\ell} \int_0^\ell \int_{-\rho}^\rho \cosh u \, du \, ds = 4\theta \ell e^{-\ell} \sinh \rho~,
\]
demonstrating the third item.
\end{proof}

Fix $T$ and let $\cR = \sqcup_{\gamma\in\cB_n(T)} R_\gamma$, whose components are \emph{return regions} of the geodesics in $\cB_n(T)$.
Disjointness of the return regions allows one to control $B_n(T)$ via $\vol\cR$: 
\begin{equation}
\label{eq:volume bound}
B_n(T) =\sum_{\gamma\in\cB_n(T)} \ell(\gamma)e^{-\ell(\gamma)} \lesssim \sum_{\gamma\in\cB_n(T)} \vol R_\gamma = \vol \cR=V_n\,\mu(\cR)~.
\end{equation}

In order to bound $\mu(\cR)$, we employ a quantification of mixing properties of the geodesic flow that follows from Selberg's celebrated $3/16$-theorem \cite{Selberg_fourier-coefficients-modular-forms_1965}. 
We reproduce the needed bounds in a convenient form for our purposes, with a short proof that relies heavily on estimates of Golubev-Kamber \cite{Golubev-Kamber_cutoff-hyperbolic-surfaces_2019} (see also \cite[Cor.~1]{Ratner_mixing-geodesic-horocycle_1987}).

\begin{proposition}
\label{thm:Selberg}
Let $Y$ be a connected finite-area hyperbolic surface whose non-zero Laplace spectrum is bounded below by $3/16$. Let $\sigma$ be normalized area measure on $Y$, and let $\mu$ be normalized Liouville measure on $T^1Y$. For real-valued $f,g\in L^2(Y,\sigma)$ and $F=f\circ \pi$ and $G=g\circ \pi$, we have 
for every $t\ge0$ that:
\[
\left| \int_{T^1Y} F(v) G(g^tv)\, d\mu(v) -\sigma(f)\sigma(g)\right| \le (1+t) e^{-t/4} \|f-\sigma(f)\|_2\|g-\sigma(g)\|_2.
\]
\end{proposition}

\begin{proof}
By \cite[Prop.~7.2]{Golubev-Kamber_cutoff-hyperbolic-surfaces_2019}, the assumption about the spectrum implies that the spherical averaging operator $A_t$ has bounded norm on the mean-zero subspace of $L_2(Y)$: 
\begin{equation}
\label{eq:operator-norm}
\|A_t\|_{L_0^2(Y)\to L_0^2(Y)}\le (1+t)e^{-t/4}~.
\end{equation}
The integral in question is the $L^2$-pairing of $f$ and $A_tg$, i.e.
\[
\int_{T^1Y} F(v) G(g^tv)\, d\mu(v) = \langle f,A_tg\rangle_{L^2(Y)}~,
\]
hence the desired inequality follows from \eqref{eq:operator-norm} by adjusting $f$ and $g$ to $f-\sigma(f),g-\sigma(g)\in L^2_0(Y)$ and applying Cauchy-Schwarz.
\end{proof}

In the proof below, we apply this proposition to indicator functions of (small) disks.

\begin{lemma}
\label{lem:B_n bound medium}
There are constants $C, T_0>0$, independent of $n$, so that, for $T\ge T_0$ we have
\begin{align}
B_n(T) & \le CV_n(1+T) e^{-T/4}+C.
\end{align}
\end{lemma}

\begin{proof}
Fix $C,r,\delta>0$ as in Lemma~\ref{lem:tubular nbhd}.

For each $x\in X$, define $f_x=\mathbbm{1}_{B_X(x,r)}$, the indicator function of $B_X(x,r)$, and $F_x=f_x\circ \pi$. Define the overlap function
\[
K_r(y,z)=\int_X f_x(y)f_x(z) \, dA(x)~.
\]
Of course, $K_r(y,z)=\area(B_X(y,r)\cap B_X(z,r))$.

Suppose that $y\in\pi(\cR)$, so that $B_X(y,r)$ embeds in $X$. If $d_X(y,z)<r/2$, then the embedded disk $B_X(y,r/4)$ lies in $B_X(y,r)\cap B_X(z,r)$. Using the formula for area of a hyperbolic disk, 
\[
K_r(y,z) \ge 2\pi(\cosh(r/4)-1)\ge \frac{\pi r^2}{16}~.
\]
Therefore for $v\in\cR_\gamma$ and $|t-\ell(\gamma)|\le \delta$ we find
\[
\int_{\ell(\gamma)-\delta}^{\ell(\gamma)+\delta} K_r(\pi v,\pi g^t v)\, dt \ge \frac{\pi r^2\delta}8~.
\]
Every such $t$ is contained in $J_T:=[T-\delta,T+1+\delta]$. Therefore, pointwise on $T^1X$,
\[
\mathbbm{1}_\cR(v) \le \frac8{\pi r^2\delta}\int_{J_T} K_r(\pi v,\pi g^t v)\, dt~.
\]
Integrating with respect to $\mu$ and interchanging integral orders we find
\[
\mu(\cR) \lesssim \int_{J_T}\int_X \int_{T^1X} F_x(v)F_x(g^tv)\, d\mu(v) \,dA(x)\, dt~.
\]
By Selberg's Theorem \cite{Selberg_fourier-coefficients-modular-forms_1965}, Proposition~\ref{thm:Selberg} applies to $X$, so there is a constant $C_\delta$ independent of $n$ such that:
\[
\left| \int_{T^1X} F_x(v)F_x(g^tv) \, d\mu(v) \right| \lesssim C_\delta (1+T)e^{-T/4} \|f_x-\sigma(f_x)\|_2^2+\sigma(f_x)^2.
\]
It is straightforward to verify that $\|f_x-\sigma(f_x)\|_2^2\le \sigma(f_x)$,
so we find that
\[
\left| \int_{T^1X} F_x(v)F_x(g^tv) \, d\mu(v) \right| 
\lesssim (1+T)e^{-T/4} \sigma(f_x) + \sigma(f_x)^2~.
\]
Integrating over $X$ and $J_T$, and using $\sigma(f_x)\lesssim 1/\area(X)$ we find
\begin{align*}
\mu(\cR) & \lesssim 
\area(X) \left( (1+T)e^{-T/4} \frac1{\area(X)} + \frac1{\area(X)^2} \right)  \\
&= (1+T)e^{-T/4} + \frac1{\area (X)}~.
\end{align*}
Because $V_n=2\pi\area(X)$, \eqref{eq:volume bound} implies that
\[
B_n(T) \lesssim V_n (1+T)e^{-T/4} + \frac{V_n}{\area(X)}\lesssim V_n (1+T)e^{-T/4} +1~.
\]

\end{proof}

Observe that \eqref{eq:B_n2}, demonstrating control over the medium regime, follows immediately from Lemma~\ref{lem:B_n bound medium} using the estimate $V_n\asymp n^3$: The function $(1+x)e^{-x/4}$ is decreasing for $x\ge 3$, and $T$ is in the medium regime, so $T\ge 16\log n$. Hence
\[
V_n(1+T)e^{-T/4} \le V_n (1+16\log n)e^{-4\log n} =\frac{V_n}{n^4}(1+16\log n)~.
\]
In particular, $V_n(1+T)e^{-T/4}=O\left(\frac{\log n}{n}\right)$, so 
\[
\sum_{T\text{ is medium} }B_n(T) \lesssim \sum_{T=\lceil 16\log n\rceil}^{\lceil \kappa n \log n \rceil} 1 = O\left(n \log n\right)~.
\]

%% file: multiplicity-paper/proofs.tex
\section{Proof of the main Theorem~\ref{thm:m lower bounds}}
\label{sec:proofs}

We show how Theorem~\ref{thm:proportion-McShane-gaps} implies Theorem~\ref{thm:m lower bounds}.
Throughout this section, $c$ is a fixed cusp of $X(n)$, and we suppress it from the notation: $h_{\ne}(n)=h_{\ne}(X(n),c)$.

\begin{proof}[Proof of Theorem~\ref{thm:m lower bounds}]
We start with the output of Theorem~\ref{thm:proportion-McShane-gaps}: for each $\delta>0$, for large enough $n$ we have
$h_{\ne}(n)>1-\delta$.

Recall \eqref{eq:McShane2} of Lemma~\ref{lem:McShane Xn}:
\[
1-\delta < h_{\ne}(n) = 2\sum_{\alpha\in\cA_{\ne}} \frac{\omega(nq_\alpha)}{nq_\alpha}~.
\]
The set of arcs $\cA_{\ne}$ can be written as a union:
\[
\cA_{\ne} =\bigsqcup_{q=1}^\infty\; \cA_{\ne}(q)~, 
\]
where $\cA_{\ne}(q)$ is the subset of $\cA_{\ne}$ of $\lambda$-length $q$.
It is immediate that $|\cA_{\ne}(q)|\le2n\mul_{\ne}(n,q)$, so
\[
\sum_{\alpha\in\cA_{\ne}} \frac{\omega(nq_\alpha)}{nq_\alpha} \le \sum_{q=1}^\infty \frac{\omega(nq)}{nq} \cdot 2n\mul_{\ne}(n,q)~.
\]
Because $\omega(x)=\frac{1}{x}+O(\frac{1}{x^3})$, for $n$ large enough we find that
\[
1-\delta< h_{\ne}(n)
< 2(2+\delta) \sum_{q=1}^\infty \frac{\mul_{\ne}(n,q)}{nq^2}~.
\]

Fix an integer $Q>2$ whose precise value will be chosen later, and let \[S_Q=\displaystyle \sup \{\mul_{\ne}(n,q) \colon q\in \bN \;, q>Q\}.\]
Using the inequalities $\mul_{\ne}(n,q)\le \frac{1}{2}\phi(q)$ for $3\le q\le Q$, and $\mul_{\ne}(n,q)\le S_Q$ for $q>Q$, we have 
\[
\frac{1-\delta}{2+\delta}<
2\sum_{q=1}^\infty \frac{\mul_{\ne}(n,q)}{nq^2} \le 
\frac1n\left(2+\frac12+\sum_{q=3}^Q  \frac{\phi(q)}{q^2} \right) + \frac{2S_Q}{n}\sum_{q=Q+1}^\infty \frac{1}{q^2}~.
\]

The growth of the first sum can be obtained from classical techniques (see \cite[\S3.10-3.12]{Apostol} and \cite[\S18.5]{Hardy-Wright_intro-theory-numbers_1938} for detail), and compare the second sum to an integral:
\[
\sum_{q=1}^Q \frac{\phi(q)}{q^2} = \frac6{\pi^2}\log Q + O(1)~, \ \ \text{ and}
\qquad
\sum_{q>Q} q^{-2} \le \int_{x>Q} \frac{dx}{x^2} = Q^{-1}.
\]
We find after multiplying by $n$ that
\[
\frac{1-\delta}{2+\delta}n<\frac6{\pi^2}\log Q+O(1) +\frac{2S_Q}Q~.
\]
Thus for every $\frac{\pi^2}{12}>\varepsilon>0$, we have 
for large enough $n$ and for every $Q$ that
\begin{equation}
\label{eq:split-at-Q}
\left(\frac{\pi^2}{12}-\frac{\varepsilon}{2}\right)n
< \log Q + \frac{\pi^2}3\frac{S_Q}{Q}~,
\end{equation}
namely
\[
S_Q>\frac {3Q}{\pi^2}\left(\left(\frac{\pi^2}{12}-\frac{\varepsilon}{2} \right) n - \log Q\right)~.
\]
We may take $Q=\left\lfloor
\exp\left( \left(\frac{\pi^2}{12}-\varepsilon\right)n\right)\right\rfloor$ to obtain:
\[
S_Q>\frac{3Q}{\pi^2}\left(\left(\frac{\pi^2}{12}-\frac{\varepsilon}2\right)n-\log Q\right) \ge \frac{3\varepsilon}{2\pi^2}nQ~.
\]
and $n$ large enough, we have $S_Q> \exp\left(\left(\frac{\pi^2}{12}-\varepsilon\right)n\right)$, as desired.

To finish, we explain why the large multiplicity can be realized for arcs with denominators that are not too long (that is, roughly exponential in $n$), by taking more care with the tail of $h_{\ne}(n)$. 
Indeed, observe that if the arc  $\alpha(p/q)$ is simple, then the continued fraction coefficients of $p/q$ are bounded by $n$.
Given a word $w=(a_1,\ldots,a_k)\in\{1,\ldots,n\}^k$ with denominator $q(w)$, Lemma~\ref{lem:continued fraction computation} implies 
\[
q(w)^{-2}\lesssim 2|I_w|~, \ \ \text{ and } \ \ 
q(w)\lesssim \prod_i(1+a_i)~.
\]
Now Lemma~\ref{lem:continued fraction harder computation} implies that there is a constant $\beta>0$ so that
\[
\sum_k \sum_{w\in\{1,\ldots,n\}^k} q(w)^{-2+\frac{\beta}{n+1}} \lesssim \sum_k \left( \frac{2n+1}{2n+2}\right)^k \lesssim n~. 
\]
Let $\cW(Q)$ be the set of words $w\in\{1,\ldots,n\}^k$, so that the denominator $q(w)$ is at least $Q$.
Each rational $p/q\in\cS\bQ(n)\cap[0,1/2]$ with $q>Q$ determines one word in $\cW(Q)$ (say, with final digit at least $2$), so we find that
\[
\sum_{q> Q} \frac{\mul_{\ne}(n,q)}{q^2} \le \sum_{w\in\cW(Q)} q(w)^{-\frac{\beta}{n+1}} q(w)^{-2+\frac{\beta}{n+1}} \le Q^{-\frac{\beta}{n+1}} \sum_{k,w\in\{1,\ldots,n\}^k}q(w)^{-2+\frac{\beta}{n+1}} \lesssim nQ^{-\frac{\beta}{n+1}}~.
\]

Now we repeat the above argument culminating in
\eqref{eq:split-at-Q}, with the slight flexibility of splitting the sum for $h_{\ne}(n)$ into three pieces. 
We find that there is another absolute constant $C>0$ so that, for every $\varepsilon>0$, for $n$ large enough, and for every $Q_1<Q_2$,
\[
\left(\frac{\pi^2}{12}-\frac\varepsilon2\right)n < \log Q_1 + \frac{\pi^2}{3} \frac{S_{Q_1,Q_2}}{Q_1} + CnQ_2^{-\beta/(n+1)}~,
\]
where $S_{Q_1,Q_2}$ is the maximum of $\mul_{\ne}(n,q)$ among $Q_1\le q\le Q_2$.
Now we choose 
\[
Q_1=\left\lfloor
\exp\left( \left(\frac{\pi^2}{12}-\varepsilon\right)n\right)\right\rfloor~, \ \ \text{ and } \ \ Q_2=\left\lfloor
\exp\left( B(n+1)\right)\right\rfloor~,
\]
for a large exponent $B$, to be chosen momentarily. The above inequality takes the form
\[
\left(\frac{\pi^2}{12}-\frac\varepsilon2\right)n < \left(\frac{\pi^2}{12}-\varepsilon\right)n + \frac{2\pi^2}{3} e^{-\left(\frac{\pi^2}{12}-\varepsilon\right)n}
S_{Q_1,Q_2} + Cne^{-\beta B}~,
\]
(note the trivial bound $Q_1^{-1}\le 2\exp(-(\frac{\pi^2}{12}-\varepsilon)n)$, and a slightly larger $C$), so that
\[
S_{Q_1,Q_2}> \frac{3}{2\pi^2}\left( \frac{\varepsilon}2 -Ce^{-\beta B} \right) n \;  e^{\left(\frac{\pi^2}{12} -\varepsilon\right)n}~.
\]
Evidently, if $B$ is chosen large enough (dependent on $\varepsilon$) 
\[
S_{Q_1,Q_2}> e^{\left(\frac{\pi^2}{12} -\varepsilon\right)n}~,
\]
as claimed.
\end{proof}

%% file: multiplicity-paper/Mirzakhani.tex
\section{Simple length multiplicity via Mirzakhani's asymptotic}

It is instructive to reflect on the following failed approach to producing large simple length multiplicity in $X(n)$ which relies on the pigeonhole principle, using Mirzakhani's asymptotic formula for the growth of the number of simple closed geodesics and the fact that (absolute values of) traces of elements of $\Gamma(n)$ are integral.

Let us first recall Mirzakhani's asymptotic formula.
Consider a complete hyperbolic surface $Y$ of genus $g\in \bN$ with $\nu\in \bN$ punctures, and let $d=6g-6+2\nu>0$.
Mirzakhani shows that the number $N_Y(L)$ of simple closed geodesics in $Y$ of length at most $L\in \bR_{\ge 0}$ is
\[
N_Y(L)\sim c(Y) \cdot L^d
\]
for some constant $c(Y)$ which one can describe explicitly \cite{Mirzakhani_Simple-Geodesics-Volumes_2007,Mirzakhani_growth-number-simple-godesics-hyperbolic-surfaces_2008}.

Now we specialize to $X(n)$, and write $d(n)=6g(n)-6+2\nu(n)$ and $c(n)$ for the constants appearing in the asymptotic.
Observe that conjugacy classes in $\Gamma(n)$ have integral traces: for $t\in \bN$, let $\mul_{\tr}(n,t)$ be the number of simple closed geodesics on $X(n)$ of trace $\pm t$ and $\mul_{\tr}(n)=\sup \{\mul_{\tr}(n,t)\colon t>2\}$. We have
\(
N_{X(n)}(L)=\sum_{3\le t\le \lfloor2\cosh(L/2)\rfloor} \mul_{\tr}(n,t)
\)
so the pigeonhole principle yields:
\begin{equation*}
\max \{\mul_{\tr}(n,\lvert \tr \rvert(\gamma)) \colon \ell_X(\gamma)\le L\}
\ge \frac{N_{X(n)}(L)}{\lfloor 2\cosh(L/2)\rfloor-2} \sim c(n) \cdot L^{d(n)}\cdot e^{-L/2}
\end{equation*}

Assume that $\log c(n)\ge -\sigma\, d(n) \log(n) $for some $\sigma\in (0,3)$ (with some work, ChatGPT-Pro could show that $\log c(n)\ge (-3\log(n)-C) d(n)$ for some $C>0$; but $\sigma<3$ seems out of reach, albeit not excluded).
In that case, we would have:
\begin{equation}
\label{eq:heuristic}
\log\left(\max \{\mul_{\tr}(n, \lvert \tr \rvert(\gamma))  \colon \ell_X(\gamma)\le L\} \right)
\ge -\sigma\, d(n) \log(n)  + d(n) \log(L) - \frac{L}{2}~.
\end{equation}
Recall that $d(n)\asymp n^3$.
If we could choose $L(n)\approx n^{\delta}$ for some $\delta\in (\sigma, 3)$ then \eqref{eq:heuristic} would reveal that $\log \mul_{\tr}(n) \gtrsim d(n)\log(n)\asymp n^3\log(n)$.
Note that this growth is much faster than the inequality $\log  \mul_{\tr}(n) \ge \left(\frac{\pi^2}{12}-\epsilon\right)n$ as provided by our main Theorem \ref{thm:sup-mul>exp}.

However, even if we had the required lower bound $\log c(n)\ge (-\sigma \log(n)-C) d(n)$ this heuristic poses a fundamental problem:
The choice of $L$ would require a strong refinement of Mirzakhani's asymptotic to precise growth estimates of $N_{X(n)}(L)$, that are currently unavailable in any explicit range whatsoever.

We believe that such questions are of independent interest.